\documentclass[12pt]{amsart}

\usepackage{amssymb,amsmath}
\usepackage{graphpap}
\usepackage{supertabular}
\usepackage{array}

 
\newtheorem{theorem}[subsection]{Theorem}
\newtheorem*{theorem*}{Theorem}

\newtheorem{proposition}[subsection]{Proposition}
\newtheorem{corollary}[subsection]{Corollary}

\newtheorem*{theorem1}{Main theorem}

\theoremstyle{definition}
\newtheorem{definition}[subsection]{Definition}

\newtheorem{example}[subsection]{Example}

\theoremstyle{remark}
\newtheorem{remark}[subsection]{Remark}

\newcommand{\mt}[1]{\operatorname{#1}}
\newcommand{\EEE}{{\mathbb E}}
\newcommand{\DDD}{{\mathbb D}}
\newcommand{\AAA}{{\mathbb A}}
\newcommand{\QQ}{{\mathbb Q}}
\newcommand{\ZZ}{{\mathbb Z}}
\newcommand{\CC}{{\mathbb C}}
\newcommand{\OO}{{\mathcal O}}
\newcommand{\RR}{{\mathbb R}}
\newcommand{\PP}{{\mathbb P}}
\newcommand{\FF}{{\mathcal F}}
\newcommand{\NN}{{\mathbb N}}
\newcommand{\FFF}{{\mathbb F}}
\newcommand{\Supp}{\mt{Supp}}
\newcommand{\Sing}{\mt{Sing}}
\newcommand{\Diff}{\mt{Diff}}
\newcommand{\ord}{\mt{ord}}
\newcommand{\wt}{\mt{wt}}
\newcommand{\D}{{\Delta}}
\newcommand{\Hom}{\mt{Hom}}
\newcommand{\Exc}{\mt{Exc}}
\newcommand{\codim}{\mt{codim}}
\newcommand{\mult}{\mt{mult}}

\newcommand{\GG}{{\Gamma}}
\newcommand{\GGG}{{\Gamma_2}}
\newcommand{\GGGG}{{\Gamma_3}}
\newcommand{\GGGGG}{{\Gamma_4}}

\newcommand{\ZZZ}{{\Delta}}
\newcommand{\ZZZZ}{{\Delta_3}}
\newcommand{\ZZZZZ}{{\Delta_4}}

\newcommand{\XX}{{\Upsilon}}
\newcommand{\XXX}{{\Upsilon_4}}

\newcommand{\YY}{{\Omega}}

\newcommand{\down}[1]{\llcorner #1 \lrcorner}
\newcommand{\up}[1]{\ulcorner #1 \urcorner}
\newcommand{\fr}[1]{\{ #1\}}

\title{Classification of three-dimensional exceptional log canonical
hypersurface singularities \large{II}}
\author{S.~A.~Kudryavtsev}

\date{}
\address{Department of Algebra, Faculty of Mathematics,
Moscow State Lomonosov University, 117234 Moscow,
Russia}
\email{kudryav@mech.math.msu.su}

\begin{document}
\begin{abstract}
In this paper the detailed classification of
three-dimensional exceptional canonical
hypersurface singularities which don't satisfy
the condition of well-formedness is given.
This result completes the
classification of
three-dimensional exceptional log canonical
hypersurface singularities started in \cite{Kud3}.
\end{abstract}
\maketitle

\section*{\bf {Introduction}}
In this paper the
classification of
three-dimensional exceptional hypersurface singularities started in \cite{Kud3}
is completed.
The main classification theorem is the following one:

\begin{theorem1}
Let $(X,0) \subset(\CC^4,0)$ be a three-dimensional exceptional canonical
$($respectively strictly log canonical$)$ hypersurface singularity
defined by a polynomial $f$.
Then there exists a biholomorphic coordinate change
$\psi \colon (\CC^4,0) \to (\CC^4_{t,z,x,y},0)$ and unique primitive vector
${\bf p}\in N_{\RR}$ such that just one of the following two possibilities holds:
\begin{enumerate}
\item The quasihomogeneous polynomial $\tilde f_{\bf p}=(f\circ \psi)_{\bf p}$
defines an exceptional canonical $($respectively strictly log canonical
and canonical outside
{\bf 0}$)$ singularity $(X_{\bf p},0) \subset(\CC^4_{t,z,x,y},0)$.
In this case ${\bf p}$-blow-up of $\CC^4$ induces purely log terminal blow-ups
$\varphi \colon (Y,E)\to (X,0)$ and
$\varphi_{\bf p} \colon (Y_{\bf p},E_{\bf p})\to (X_{\bf p},0)$, where
$(E,\Diff_E(0))=(E_{\bf p},\Diff_{E_{\bf p}}(0))$.
That is, these singularities have the same type and in particular the same
complement index.
\par
The canonical singularities satisfying the condition of well-formedness ---
$\Diff_{E/\PP(\bf p)}(0)=0$ are classified in the theorems
3.24, 3.28 of \cite{Kud3} and in the tables of chapter \cite[\S 4]{Kud3}.
The polynomial
$\tilde f_{\bf p}$; $(E,\Diff_E(0))$; minimal complement index are written
in the tables.
\par
The canonical singularities not satisfying the condition of
well-formedness are classified in the tables of chapter \S 3.
The polynomial
$\tilde f_{\bf p}$; $(E,\Diff_E(0))$; minimal complement index are written
in the tables.
\par
The strictly log canonical and canonical outside
{\bf 0} quasihomogeneous singularities are always exceptional
$($in any dimension$)$ by theorem
\cite[2.13]{Kud3}.
In the three-dimensional case
$(E,\Diff_E(0))=(\tilde f_{\bf p}\subset\PP({\bf p}),0)$ is
$K3$ surface with Du Val singularities and $(X,0)$ is
1-complementary.
\item $\tilde f_{\bf p}=t^3+g_2^2(z,x,y)$, where $g_2$ is an irreducible
homogeneous polynomial of degree two.
In this case the purely log terminal blow-ups are constructed in the theorems
\cite[3.3, 3.5]{Kud3}. Also it was obtained the similar
classification depending on the type of jets
$\tilde f_5$ and $\tilde f_6$.
\end{enumerate}

\begin{proof}
According to the main theorem proved in
\cite{Kud3} it is enough to classify the three-dimensional exceptional
canonical quasihomogeneous singularities which are not well-formed.
The required classification is given in the tables of chapter \S 3.
\end{proof}
\end{theorem1}

It follows from the classification that there are only finite number of
types of three-dimensional exceptional log canonical
hypersurface singularities.
There is a conjecture that the similar result about
the finite number of types is true in any dimension for the
exceptional singularities.
See also examples \cite[1.19, 1.20]{Kud3}.
\par
The present paper is the direct continuation of
\cite{Kud3}. Therefore see \cite{Kud3} in connection with the notations,
definitions and etc.
\par
The paper is organized in the following way. In chapter
\S 1 the purely log terminal blow-ups of not well-formed singularities
are studied.
In chapter \S 2 the three-dimensional quasihomogeneous not well-formed
singularities are considered.
In chapter \S 3 the summarizing tables are written.

I am grateful to Professor Yu.G.Prokhorov for
useful discussions and valuable remarks.
\par
The research was partially
supported by a grant 99-01-01132 from the Russian Foundation of Basic Research
and a grant INTAS-OPEN 2000\#269.
\par

\section{\bf Structure of purely log terminal blow-up of not
well-formed singularity}
Now we explain in detail the condition of well-formedness of quasihomogeneous
singularity.
\par
Let $(X,0)\subset (\CC^n,0)$ be a quasihomogeneous log canonical,
canonical outside {\bf 0}
singularity with the integer weights
${\bf p}=(p_1,\ldots,p_n)$ and quasihomogeneous degree $d$ defined by a
polynomial $g(x_1,\ldots,x_n)$. Without loss of generality it can be assumed that
$(p_1,\ldots,p_n)=1$.
\par
Consider ${\bf p}$-blow-up $\psi\colon \CC^n({\bf p})\to \CC^n$.
The exceptional divisor of $\psi$ is the weighted projective space
$\PP({\bf p})$. This blow-up induces a plt blow-up
$\psi|_{X{\bf p}}\colon (X({\bf p}),E)\to (X,0)$ (see theorem 2.13 \cite{Kud3}).
Put $q_i=(p_1,\ldots,\hat{p}_i,\ldots,p_n)$. Then
$E=(\tilde g(x_1,\ldots,x_n)\subset \PP(\tilde p_1,\ldots,\tilde p_n))=
(g(x_1^{1/q_1},\ldots,x_n^{1/q_n})\subset
\PP(p_1q_1/(q_1\cdots q_n),\ldots,p_nq_n/(q_1\cdots q_n))$. Put
$q_{ij}=(\tilde p_1,\ldots,\widehat{\tilde{p}}_i,\ldots,
\widehat{\tilde{p}}_j,\ldots,\tilde p_n)$
if $i\ne j$.
The quasihomogeneous degree of $\tilde g(x_1,\ldots,x_n)$
is denoted by $\tilde d$.

\begin{definition}\label{w-d}\cite{Fletcher}
In the above--mentioned notations the quasihomogeneous singularity is called
{\it well--formed singularity} if one of the following equivalent conditions
holds:
\begin{enumerate}
\item $\Diff_{E/\PP({\bf p})}(0)=0$;
\item One have $q_{ij} | \tilde d$ for all $i\ne j$.
\end{enumerate}
In opposite case our singularity is called
{\it not well--formed singularity}. Let us remark that if we choose the other
quasihomogeneous weights then it can happen that the condition
of well-formedness is not
fulfilled. For the weakly exceptional and consequently for the exceptional
singularity the primitive vector
${\bf p}$ is unique \cite[2.13]{Kud3} and therefore the notion of
well-formedness is uniquely defined.
\end{definition}

Put
$I_i=\big\{\ j| 1\le j\le n,\ j\ne i,\ q_{ij}\nmid \tilde d\ \big\}.$

\begin{proposition}\label{diff}
$$\Diff_E(0)=\sum_{i=1}^n\Big(1-\frac1{q_i}\Big)
\Big(C_i\setminus\bigcup_{j\in I_i}C_{ij}\Big)+
\sum_{i<j; \ j\in I_i}
\Big(1-\frac1{q_{ij}\cdot q_i\cdot q_j}\Big) C_{ij},$$
where the divisors $C_i$ and $C_{ij}$ on $E$ are defined by the equations $x_i=0$
and $x_i=x_j=0$ respectively. Moreover
$C_i'\stackrel{\rm def}{=}C_i\setminus\bigcup_{j\in I_i}C_{ij}$=
$\big\{x_i=\big(\prod_{j\in I_i}x_j\big)^{-1}\cdot
g(x_1,\ldots,x_{i-1},0,x_{i+1},\ldots,x_n)=0\big\}$, where $q_i\ne 1$.

\begin{proof} To calculate the different we have to determine the codimension
2 components of
$\Sing X({\bf p})$ in $X({\bf p})$ which are lying in $E$.
The variety $X({\bf p})$ in the affine piece
$U_k=\CC^n_{x_1,\ldots,x_n}\big/
\ZZ_{p_k}(-p_1,\ldots,-p_{k-1},1,-p_{k+1},\ldots,-p_n)$ of
${\bf p}$-blow-up is given by a polynomial
$g'=g(x_1,\ldots,x_{k-1},1,x_{k+1},\ldots,x_n)=0$.
It is clear that the singularity defined by this polynomial $g'$ has
the codimension not less then 3 in
$\CC^{n-1}\cong\{x_k=0\}$ and hence doesn't influence on
$\Diff_E(0)$ (recall that $E\bigcap U_k$
is given by a polynomial $g'$ in $\{x_k=0\}/\ZZ_{p_k}$). Thus $\Diff_E(0)$
appears only because of the action of cyclic group.
Since
$q_{ij}\cdot q_i\cdot q_j=(p_1,\ldots,\hat p_i,\ldots,
\hat p_j,\ldots, p_n)$ then it remains to show the following statement: if
$q_{ij}\nmid \tilde d$ then $\{x_i=x_j=0\}$ defines a divisor on $E$. Actually,
this condition means that
$x_i$ or $x_j$ divides any monomial of
$g$.
\end{proof}
\end{proposition}

\begin{remark}\label{view}
If $q_{ij}\nmid \tilde d$ then
$g=x_ig_1(x_1,\ldots,x_n)+x_jg_2(x_1,\ldots,x_n)$.
Moreover, if $q_i>1$ and $q_{ij}\nmid \tilde d$ then $q_j=1$. Indeed,
if $q_j>1$ then
$g=x_i^{q_i}g_1(x_1,\ldots,x_n)+x_j^{q_j}g_2(x_1,\ldots,x_n)$. Hence $g$
defines a nonnormal singularity.
\end{remark}

The next proposition is proved in the same way as the previous proposition
\ref{diff}.
\begin{proposition}\label{diffp}
$$\Diff_{E/\PP({\bf p})}(0)=
\sum_{i<j; \ j\in I_i}
\Big(1-\frac1{q_{ij}}\Big) C_{ij}.$$
\end{proposition}

\begin{corollary}$($cf. \cite[2.10]{I}$)$ Let $(X,0)\subset (\CC^n,0)$ be a
quasihomogeneous strictly log canonical, canonical outside
{\bf 0} singularity with weight
{\bf p}. Then it is well-formed (by theorem 2.13 \cite{Kud3}
it is exceptional).
\begin{proof} Let
$\psi\colon (X({\bf p}),E)\to (X,0)$ be a plt blow-up.
Since $K_{X({\bf p})}+E=\psi^*K_X$ is Cartier
divisor then $\Diff_E(0)=0$. Hence by propositions \ref{diff}
and \ref{diffp}
$\Diff_{E/\PP({\bf p})}(0)=0$. In particular $E$ is
$K3$ surface with Du Val singularities if
$n=4$ \cite[1.4.1]{Fletcher}.
\end{proof}
\end{corollary}

\begin{example}\label{ex1}
Let $(X,0)\subset (\CC^3,0)$ be a weakly exceptional Du Val singularity
defined by a polynomial
$h$. Then $h$ is one of the following polynomials:
$x_1^2+x_2^2x_3+x_3^{n-1}$,\ $n\ge 4$ (type $\DDD_n$);
$x_1^2+x_2^3+x_3^4$ (type $\EEE_6$); $x_1^2+x_2^3+x_3^3x_2$ (type $\EEE_7$);
$x_1^2+x_2^3+x_3^5$ (type $\EEE_8$) (see \cite[4.7]{Pr2}).
Consider ${\bf p}$-blow-ups corresponding to their quasihomogeneous weights.
Then the singularities $\DDD_n$, $n\ge 5$ and $\EEE_7$ are not well-formed.\\
\begin{center}
\renewcommand{\arraystretch}{1.5}
\begin{tabular}{|c|c|c|c|}
\hline
Singularity &$E\cong \PP^1$&$\Diff_{E/\PP({\bf p})}(0)$&$\Diff_E(0)$\\
\hline
$D_{2n},\ n\ge 3$& $x_1+x_2^2x_3+x_3^{2n-1}\subset $&
$\frac{n-2}{n-1}P_3$&$\frac12P_1+\frac12P_2+$\\
&$\subset\PP(2n-1,n-1,1)$&&$+\frac{2n-3}{2n-2}P_3$\\
\hline
$D_{2n+1},\ n\ge 2$& $x_1^2+x_2x_3+x_3^{2n}\subset $&
$\frac{2n-2}{2n-1}P_3$&$\frac12P_1+\frac12P_2+$\\
&$\subset\PP(n,2n-1,1)$&&$+\frac{2n-2}{2n-1}P_3$\\
\hline
$E_7$& $x_1+x_2^3+x_3x_2\subset $&
$\frac12P_3$&$\frac12P_1+\frac23P_2+$\\
&$\subset\PP(3,1,2)$&&$+\frac34P_3$\\
\hline
\end{tabular}
\end{center}
\end{example}

\begin{definition}\label{def1}
Let us define a $\QQ$-divisor $\widehat D=\sum^n_{i=1}(1-\frac1{q_i})\{x_i=0\}$
on
$\PP({\bf p})$ and a $\QQ$-divisor $D=\sum^n_{i=1}(1-\frac1{q_i})C_i$ on $E$.
Recall that, if $(X,0)$ is a well-formed singularity then
$D=\widehat D|_E=\Diff_E(0)$ (cf. \cite[3.12]{Kud3}).
Nevertheless the following proposition takes place in the general case.
\end{definition}

\begin{proposition} $K_E+\Diff_E(0)=(K_{\PP({\bf p})}+E+\widehat D)|_E$.
\begin{proof} By proposition \ref{diffp}, remark \ref{view}
and corollary \cite[3.10]{Sh1}
$(K_{\PP({\bf p})}+E+\widehat D)|_E=K_E+\Diff_{E/\PP({\bf p})}(\widehat D)=
K_E+\sum_{i<j; \ j\in I_i}\Big[(1-\frac1{q_{ij}})C_{ij}+
\frac1{q_{ij}}\cdot(1-\frac1{q_i\cdot q_j})C_{ij}\Big]+
\sum^n_{i=1}(1-\frac1{q_i})(C_i\setminus\bigcup_{j\in I_i}C_{ij})$.
By proposition \ref{diff} the last expression is equal to $K_E+\Diff_E(0)$.
\end{proof}
\end{proposition}

\begin{example} Let us return to the example \ref{ex1}. Then
\begin{center}
\renewcommand{\arraystretch}{1.5}
\begin{tabular}{|c|c|c|}
\hline
Singularity& $D$ & $\Diff_E(0)=\Diff_{E/\PP({\bf p})}(\widehat D)$\\
\hline
$D_{2n},\ n\ge 3$& $\frac12P_1+\frac12P_2+\frac12P_3$
&$\frac12P_1+\frac12P_2+$\\
& &$+\big(\frac{n-2}{n-1}+\frac1{n-1}\cdot\frac12\big)P_3$\\
\hline
$D_{2n+1},\ n\ge 2$& $\frac12P_1+\frac12P_2$
&$\frac12P_1+\frac12P_2+$\\
&&$+\big(\frac{2n-2}{2n-1}+\frac1{2n-1}\cdot 0\big)P_3$\\
\hline
$E_7$& $\frac12P_1+\frac23P_2+\frac12P_3$
&$\frac12P_1+\frac23P_2+$\\
&&$+\big(\frac12+\frac12\cdot\frac12\big)P_3$\\
\hline
\end{tabular}
\end{center}
\end{example}

\section{\bf Investigation of three-dimensional not well-formed singularities
on exceptionality}

The next theorem is proved by exhaustion of all cases. This sorting is described
in \cite{Kud3}.
\begin{theorem} All three-dimensional not well-formed exceptional canonical
quasihomogeneous singularities can be obtained by the rotations of Newton's line
passing through the following parts of type $\mathcal M_2$:
\par $(1)$
$t^2+z^3+x^m$, where
$x^m$ --- $zx^5$, $zx^5y$, $zx^7$, $zx^6y$, $zx^5y^2$.
\par $(2)$
$t^2+z^4+x^m$, where
$x^m$ ---  $zx^4$, $zx^5$, $zx^4y$.
\par $(3)$
all parts of type $\Upsilon^{(2)}_2$ $($they have the equation of the form
$t^2+z^3x+h(z,x,y))$.
\par $(4)$
$t^2+z^3y+x^m$, where
$x^m$ --- $x^7y$, $x^8y$, $zx^5$, $zx^6$, $zx^5y$.
\par $(5)$
$t^2+z^5+zx^5$; $t^2+z^4x+x^m$, where
$x^m$ ---
$x^6$, $x^5y$, $x^4y^2$,
$zx^5$, $zx^4y$, $zx^3y^2$;
$t^2+z^4y+x^m$, where
$x^m$ ---
$x^5y$,
$zx^5$, $zx^4y$;
$t^2+z^3x^2+x^m$, where
$x^m$ ---
$zx^5$, $zx^4y$,
$zx^3y^2$;
$t^2+z^3y^2+zx^5$;
$t^2+z^3xy+x^m$, where
$x^m$ ---
$x^6$,
$x^5y$, $x^7$, $x^6y$, $x^5y^2$,
$zx^5$, $zx^4y$;
$t^2+z^5+zx^3y$; $t^2+z^4x+zx^3y$;
$t^2+z^4x+x^4y$; $t^2+z^3x^2+zx^3y$.
\par $(6)$
$t^3+z^2xy+x^5$; $t^3+z^2xy+x^4y$; $t^3+z^4+zx^2y$; $t^3+z^4+tx^3$;
$t^3+z^3x+x^m$, where
$x^m$ ---
$zx^2y$, $x^3y$, $tx^3$, $tx^2y$;
$t^3+z^3y+tx^3$;
$t^3+tz^3+x^m$, where
$x^m$ ---
$z^2x^2$, $zx^3$, $zx^2y$;
$t^3+tz^2x+zx^2y$; $t^3+tz^2x+x^2y^2$; $g_3(t,z)+zx^3$, where $g_3$ is a binary
form of degree 3;
$t^3+z^2x+x^m$, where
$x^m$ ---
$x^4$, $x^5$, $x^4y$,
$x^6+tx^4$, $x^5y$,
$x^4y^2$, $x^7$, $x^6y$, $x^5y^2$,
$x^4y^3$, $x^8$, $x^7y$, $x^6y^2$,
$x^5y^3$,
$x^4y^4$, $tx^3$, $tx^3y$, $tx^5$,
$tx^4y$, $tx^3y^2$; $t^3+z^2y+x^7y$; $t^3+z^2y+tx^5$.
\par $(7)$
$t^2z+g_4(z,x)$, where $g_4$ is a binary form of degree 4; $t^2z+z^4+x^m$, where
$x^m$ ---
$x^3y$, $x^5$, $x^4y$, $x^3y^2$, $tx^3$,
$tx^2y$;$t^2z+z^3x+x^m$, where
$x^m$ ---
$x^3y$, $x^5$, $x^4y$, $x^3y^2$,
$tx^3$, $tx^2y$;
$t^2z+z^3y+x^m$, where
$x^m$ ---
$x^5$, $x^4y$, $zx^3$,
$z^2x^2$, $tx^3$;
$t^2z+z^2x^2+x^m$, where
$x^m$ ---
$x^5$, $x^4y$, $x^3y^2$,
$tx^2y$; $t^2z+z^2xy+x^5$; $t^2z+z^2xy+x^4y$;
$t^2z+z^2x+x^m$, where
$x^m$ ---
$x^4$, $zx^3+x^5$, $x^4y$, $x^6$, $x^5y$, $x^4y^2$, $tx^3$,
$tx^4$, $tx^3y$;
$t^2z+z^2y+x^m$, where
$x^m$ ---
$x^5y$, $zx^3$, $tx^4$.
\par $(8)$
$t^2x+z^4+x^m$, where
$x^m$ ---
$x^3y$, $x^5$, $x^4y$, $zx^2y$, $zx^4$,
$zx^3y$, $zx^2y^2$;
$t^2x+z^3x+x^m$, where
$x^m$ ---
$x^2y^2$, $x^5$, $x^4y$,
$x^3y^2$, $x^2y^3$, $x^6$, $x^5y$, $x^4y^2$,
$x^3y^3$, $x^2y^4$, $zx^2y$, $zx^4$, $zx^3y$,
$zx^2y^2$;
$t^2x+z^3y+x^m$, where
$x^m$ ---
$x^5$, $x^4y$, $x^6$,
$x^5y$, $x^4y^2$, $zx^4$, $zx^3y$, $z^2x^2$.
\par $(9)$
$t^2y+z^4+zx^4$;
$t^2y+z^3x+x^m$, where
$x^m$ ---
$x^5$, $x^6$, $x^5y$, $zx^4$, $tx^3$.
\end{theorem}

In investigating given singularities on the exceptionality it is convenient to
use the proposition
\cite[3.13]{Kud3} and corollary \cite[3.14]{Kud3}. In general case corollary
\cite[3.15]{Kud3} is formulated in the following way.

\begin{proposition}\label{DelPezzo3}
Let $(X,0)\subset (\CC^4_{x_1,\ldots,x_4},0)$ be a three-dimensional
quasihomogeneous log canonical singularity with the weight
{\bf p}. Assume that $(X,H)$ is not lc for any hyperplane section
$H$.
Let $D=\sum d_iD_i$ (see definition \ref{def1}).
If $d_k\ge \frac67$ for some
$k$ then $(X,0)$ is exceptional.
\begin{proof} Suppose that $(X,0)$ is not exceptional.
Let $H_X=\{x_k=0\}|_X$ be a hyperplane section and let
$H_{X({\bf p})}$ be a proper transform of $H_X$ on $X({\bf p})$.
Then $K_{X({\bf p})}+E+H_{X({\bf p})}=\psi^*(K_X+H_X)+aE$.
Let $C_k=\sum_i C_{ki}$, where $C_{ki}$ are the irreducible curves
(see proposition \ref{diff}).
One have $\Diff_E(0)=\sum_i d_{ki}C_{ki}+\Xi$. The coefficients $d_{ki}\ge 6/7$ for all $i$
since
$d_{ki}=\frac{m_i-1}{m_i}+\frac1{m_i}\cdot d_k\ge d_k$ \cite[3.10]{Sh1}.
By proposition \cite[3.14]{Kud3} there exists 1,2,3,4 or 6-complement $D^+$.
By the definition of complement we have $D^+\ge \sum_i C_{ki}+\Xi$. Therefore
$-(K_E+\sum_i C_{ki}+\Xi)$ is nef and lc.
Since $\Diff_E(H_{X({\bf p})})=\Diff_{E/\PP({\bf p})}
\big(\widehat D+\{x_k=0\}\big)$, where $\{x_k=0\}$ is the corresponding divisor
on $\PP({\bf p})$ then $\Diff_E(H_{X({\bf p})})=\sum_i C_{ki}+\Xi$.
Hence $a\ge 0$ and $(X,H_X)$ is lc. This contradiction concludes the proof.
\end{proof}
\end{proposition}

\begin{remark} Proposition \ref{DelPezzo3} is hypothetically true in any
dimension
(see \cite[3.16]{Kud3}). In the two-dimensional case we must require that
$d_k\ge \frac23$ for some $k$.
\end{remark}

\begin{example} \label{exam1}
In proposition \ref{DelPezzo3} we can't change
$D$ on $\Diff_E(0)$. Let us show it in the following two
examples.
\par (1) Let $(X,0)\subset (\CC^3,0)$ be a Du Val singularity of type
$\DDD_n$ defined by a polynomial $x_1^2+x_2^2x_3+x^{n-1}_3$, where $n\ge 4$.
The pair $(X,H)$ is not lc for any hyperplane section $H$. If
$n\ge 5$ then one of the coefficients of $\Diff_E(0)$ is not less then $2/3$ and
all coefficients of $D$ are equal to $1/2$.
This singularity is not exceptional by the definition
since $K_X+\{x_1=x_3=0\}\sim_{\QQ} 0$ is 2-complement
(see also example \cite[1.8]{Kud3}).
\par (2) Let $(X,0)\subset (\CC^4_{t,z,x,y},0)$ be a canonical singularity
defined by a polynomial $t^3+z^2x+x^4+xy^5$.
The pair $(X,H)$ is not lc for any hyperplane section $H$.
The variety $X$ has $cA_2$ singularity along the curve
$C=\{t=x=z^2+y^5=0\}=\Sing X$ (outside {\bf 0}).
The divisor $T=\{t=x=0\}$ on $X$ is $\QQ$-Cartier since
$3T\sim \{x=0\}$. Therefore the divisor
$(K_X+T)|_T=K_{\CC^2_{z,y}}+\frac23C$ is klt.
Thus $(X,0)$ is not weakly exceptional singularity
\cite[4.6]{Pr2}.
\par
Now consider $(40,45,30,18)$-blow-up of $\CC^4$. It induces a plt blow-up
of $(X,0)$. Then $(E,\Diff_E(0))=
\big(t+zx+x^4+xy\subset\PP(4,3,1,3),\frac23(\{t=0\}\backslash\{t=x=0\})
+\frac12\{z=0\}+\frac45\{y=0\}+\frac89\{t=x=0\}\big)=
\big(\PP^2,\frac23l_1+\frac12l_2+\frac45l_3+\frac89l_4\big)$, where $l_i$ are
the straight lines in general position.
The divisor
$K_{\PP^2}+\frac23l_1+\frac12l_2+\frac56l_3+l_4\sim_{\QQ} 0$ is klt
6-complement of minimal index. Since the divisor
$K_{\PP^2}+\Diff_E(0)+\frac19l_4$ is anti-ample and non-klt
then by criterion \cite[2.1]{Kud2} $(X,0)$ is not weakly exceptional singularity
again.
In this example all coefficients of $D$ are not more then
$\frac45$.
\end{example}

It is very interesting to compare the following example
\ref{exam2} with the previous example
\ref{exam1}(2).
\begin{example}\label{exam2}
Let $(X,0)\subset (\CC^4_{t,z,x,y},0)$ be a canonical singularity defined by
a polynomial $t^3+z^2x+tx^3+ty^5$.
The variety $X$ has $cA_1$ singularity along the curve
$C=\{t=z=x^3+y^5=0\}=\Sing X$ (outside {\bf 0}).
Consider the divisor $T=\{t=z=0\}$
on $X$. If $T$ had been $\QQ$-Cartier divisor then
the divisor $(K_X+T)|_T=K_{\CC^2_{x,y}}+\frac12C$ would have been klt and
consequently $(X,0)$ would have been not weakly exceptional singularity.
Now we prove that
$T$ is not $\QQ$-Cartier divisor and $(X,0)$ is an exceptional singularity.
\par
Consider ${\bf p}=(30,35,20,12)$-blow-up of $\CC^4$. It induces a plt blow-up
$\psi\colon (X({\bf p}),E)\to(X,0)$.
Then $E=\big(t^3+zx+tx^3+ty\subset\PP(3,7,2,6)\big)$,
$\Diff_E(0)=\frac12\{z=t^2+x^3+y=0\}+
\frac34\{t=z=0\}+\frac45\{y=0\}=
\frac12\Delta+\frac34\GGG+\frac45L$.
In the next picture we illustrate the situation of these curves and
the singularities of surface $E$.\\

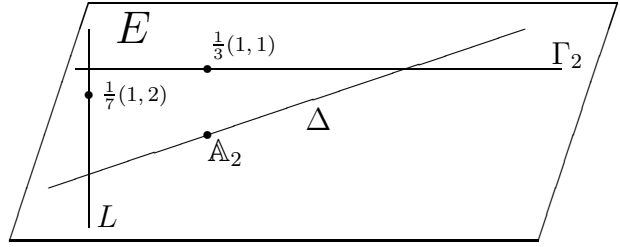
\begin{equation*}
\begin{picture}(260,70)(0,0)
\put(0,0){\line(1,3){30}}
\put(0,0){\line(1,0){200}}
\put(200,0){\line(1,3){30}}
\put(30,90){\line(1,0){200}}
\put(76,72){\tiny{$\frac13(1,1)$}}
\put(75,65){\circle*{3}}
\put(40,75){\Large{$E$}}
\put(25,65){\line(1,0){184}}
\put(205,67){$\Gamma_2$}
\put(30,5){\line(0,1){75}}
\put(30,55){\circle*{3}}
\put(35,53){\tiny{$\frac17(1,2)$}}
\put(33,5){$L$}
\put(15,20){\line(3,1){180}}
\put(75,40){\circle*{3}}
\put(75,30){$\AAA_2$}
\put(112,42){$\Delta$}
\end{picture}
\end{equation*}

Assume that $T$ is $\QQ$-Cartier divisor. Then
$K_{X({\bf p})}+E+T_Y=\psi^*(K_X+T)+aE$ is anti-ample divisor over $X$  since
$a>0$.
Hence $K_E+\Diff_E(T_Y)=K_E+\frac12\Delta+\GGG+\frac45L$ is anti-ample divisor.
It is not true since $(K_E+\Diff_E(T_Y))\cdot\Delta=
(K_E+\frac12\Delta+\frac34\GGG+\frac{11}{12}L)\cdot\Delta+(\frac14\GGG+
(\frac45-\frac{11}{12})L)\cdot\Delta=(K_{\PP({\bf p})}+E+\frac12\{z=0\}+
\frac{11}{12}\{y=0\})\cdot\Delta+\frac14-\frac7{60}=0+\frac14-\frac7{60 }>0$.
Here $\{z=0\}$ and $\{y=0\}$ are denoted the corresponding divisors on
$\PP({\bf p})$.
\par
The divisor $K_E+\Diff_{E/\PP({\bf p})}\big(\frac35\{z=0\}+\frac45\{y=0\}\big)=
K_E+\frac35\Delta+\frac45\GGG+\frac45L$ is 5-complement of minimal index.
It is easy to show that there is only the following
6-complement
$K_E+\Diff_{E/\PP({\bf p})}\big(\frac16\{t=0\}+\frac12\{z=0\}+
\frac56\{y=0\}\big)=
K_E+\frac12\Delta+\frac56\GGG+\frac56L+\frac16\{t=x=0\}\sim_{\QQ} 0$ among
1,2,3,4 and 6-complements. It is clear that this divisor is klt. Therefore
$(X,0)$ is an exceptional singularity.
\end{example}

Finally let us consider one more example showing the investigation process
of given singularities on the exceptionality.
\begin{example}
Let $(X,0)\subset (\CC^4_{t,z,x,y},0)$ be a canonical singularity defined by
a polynomial $f=t^2x+z^3x+zx^2y+atz^2y+bz^2y^3+cxy^4$.
Consider ${\bf p}=(6,4,5,3)$-blow-up of $\CC^4$.
Then $(E,\Diff_E(0))=(f\subset \PP(6,4,5,3),\frac23\Delta_3+\frac12\Upsilon_4)$,
where $\Delta_3=\{z=x=0\}$ and $\Upsilon_4=\{x=y=0\}$.
Let us prove that the divisor
$(K_{\PP({\bf p})}+E+\frac14\{z=0\})|_E=
K_E+\frac34\Delta_3+\frac12\Upsilon_4+\frac12\{t=z=0\}$ is non-klt
4-complement if
$c=0$ (then $b\ne 0$). It implies that $(X,0)$ is not exceptional singularity
if $c=0$.
Since the klt and lc properties of pair remain true by the finite dominant
morphisms then it is enough to prove that
$(Z,\frac14T)$ is lc, but not klt
pair, where $Z=\{t^2x+z^3x+zx^2+atz^2+z^2=0\}\subset
\CC^3_{t,z,x}$ and $T=\{z=0\}|_Z$. Note that $(Z,0)$ is the singularity of type
$\DDD_5$. The blow-up of $\CC^3$ with the weights $(3,4,2)$ induces a blow-up
$\psi\colon (\widetilde Z,C)\to Z$. Besides $K_{\widetilde Z}+C+
\frac14T_{\widetilde Z}=\psi^*(K_Z+\frac14T)$. The exceptional divisor $C$ is
not well-formed hypersurface $tx+zx^2+z^2\subset \PP(3,2,1)$.
The divisor $K_C+\Diff_C(\frac14T_{\widetilde Z})=K_C+\frac12P_1+(\frac12+\frac12
\cdot\frac12)P_2+(\frac23+\frac13\cdot\frac14)P_3$ is klt. By Inversion of
Adjunction we complete the proof. If $c\ne 0$ then by the same argument one can
prove that
$(X,0)$ is an exceptional singularity.
\end{example}

\section{\bf Classification of three-dimensional exceptional not well-formed
singularities}
{\bf The comments to the tables.}
All agreements of \S4 \cite{Kud3} about the notations and records in the
tables remain valid. The underlying monomials
(see point (I) of \S4 \cite{Kud3}) are written without the coefficients before
them.

\par
If log Del Pezzo surface
$E=(\widetilde f\subset \PP(\widetilde p_1,\ldots,\widetilde p_4)$ is a linear
cone, i.e. $\widetilde d=\widetilde p_k$ for some $k$ then
we write the simplified form of $E$
$$ E=\PP\Big(\frac{\widetilde p_1}{s_1},\ldots,
\frac{\widetilde p_{k-1}}{s_{k-1}},
\frac{\widetilde p_{k+1}}{s_{k+1}},\ldots,\frac{\widetilde p_4}{s_4}\Big),\
\text{where}\
s_m=\prod_{\substack{q_{ij}\nmid\widetilde d \\ i< j;\ i,j\ne m}}q_{ij}
$$

\par
For the record of different $\Diff_E(0)$ one use the following notations.
The curves $C_i'$,
where $i=1,\ldots,4$ (see proposition \ref{diff}) are denoted by
$\Gamma$, $\Delta$, $\Upsilon$, $\Omega$ respectively. The curves
$\{t=z=0\}$, $\{t=x=0\}$, $\{t=y=0\}$, $\{z=x=0\}$, $\{z=y=0\}$, $\{x=y=0\}$
are denoted by $\Gamma_2$, $\Gamma_3$, $\Gamma_4$, $\Delta_3$, $\Delta_4$,
$\Upsilon_4$ respectively.
\par
If there is a symbol * in the description of different cases then this marked
case has the other different. Then the different differs by the addition of
other curve marked by * or the different is completely another one
marked by * also.
\newpage

\begin{center}
\small{1. Singularity -- $t^2+z^3+g(z,x,y)$}
\end{center}

\begin{center}

\tabletail{\hline}
\tablehead{\hline}
\footnotesize{
\begin{supertabular}{|c|c|c|c|}
\hline
No. & $g(z,x,y)$ & Log Del Pezzo & Index \\
\hline
1& $zx^if_{5-i}(x,y^3)$ &
$t + z^3 +zx^if_{5-i}(x,y^3) \subset \PP(15,5,2,2)$  & 4\\
& & $\PP(5,1,1)$,\
$\Diff$=\scriptsize{$(0,0,0,\frac23)$}+
$\frac12\GG+\frac34\GGG$ & \\
\hline
2& $zx^5+zy^n$ &
$t + z^3 +zx+zy \subset \PP(3,1,2,2)$  & 8,10,12\\
&$n=7;9;11;13;17,19$ & $\PP^2$,\
$\Diff$=\scriptsize{$(0,0,\frac45,\frac{n-1}{n})$}+
$\frac12\GG+\frac34\GGG$ & 16,20\\
\hline
3& $zx^5+zxy^n$ &
$t + z^3 +zx^5+zxy \subset \PP(15,5,2,8)$  & 6,8,12\\
&$n=5;7;9,11;13,15$ & $\PP(5,1,4)$,\
$\Diff$=\scriptsize{$(0,0,0,\frac{n-1}{n})$}+
$\frac12\GG+\frac34\GGG$ & 16\\
\hline
4& $zx^5+azx^3y^n+bzxy^{2n}$ &
$t + z^3 +zx^5+azx^3y+bzxy^2 \subset \PP(15,5,2,4)$  & 4,8\\
&$n=3\ b\ne 0;5,7$ & $\PP(5,1,2)$,\
$\Diff$=\scriptsize{$(0,0,0,\frac{n-1}{n})$}+
$\frac12\GG+\frac34\GGG$ &\\
\hline
5& $zx^5+zx^2y^n$ &
$t + z^3 +zx^5+zx^2y \subset \PP(15,5,2,6)$  & 8,12\\
&$n=5,7;11$ & $\PP(5,1,3)$,\
$\Diff$=\scriptsize{$(0,0,0,\frac{n-1}{n})$}+
$\frac12\GG+\frac34\GGG$ &\\
\hline
6& $zx^5y+zy^n$ &
$t + z^3 +zxy+zy^n \subset \PP(3n,n,2n-2,2)$  & 8,12,16\\
&$n=7;9;13;15$ & $\PP(n,n-1,1)$,\
$\Diff$=\scriptsize{$(0,0,0,\frac{n-1}{n})$}+
$\frac12\GG+\frac34\GGG$ & 20\\
\hline
7& $zyx^if_{5-i}(x,y^2)$ &
$t + z^3 +zx^if_{5-i}(x,y^2) \subset \PP(33,11,4,2)$  & 4\\
& & $\PP(11,2,1)$,\
$\Diff$=
$\frac12\GG+\frac34\GGG$ & \\
\hline
8& $zx^5y+zxy^n$ &
$t + z^3 +g(z,x,y) \subset \PP(15n-3,5n-1,2n-2,8)$  & 8,12,16\\
&$n=6;8,10;12$ & $\PP(5n-1,n-1,4)$,\
$\Diff$=
$\frac12\GG+\frac34\GGG$ & \\
\hline
9& $zx^5y+azx^3y^{n+1}+bzxy^{2n+1}$ &
$t + z^3 +g(z,x,y) \subset \PP(15n+6,5n+2,2n,4)$  & 4,8\\
&$n=3\ b\ne 0;5$ & $\PP(5n+2,n,2)$,\
$\Diff$=
$\frac12\GG+\frac34\GGG$ &\\
\hline
10& $zx^5y+zx^2y^n$ &
$t + z^3 +g(z,x,y) \subset \PP(15n-6,5n-2,2n-2,6)$  & 8,12\\
&$n=5;9$ & $\PP(5n-2,n-1,3)$,\
$\Diff$=
$\frac12\GG+\frac34\GGG$ & \\
\hline
11& $zx^if_{7-i}(x,y)$ &
$t + z^3 +zx^if_{7-i}(x,y) \subset \PP(21,7,2,2)$  & 4\\
&$i\le 2$ & $\PP(7,1,1)$,\
$\Diff$=
$\frac12\GG+\frac34\GGG$ & \\
\hline
12& $zx^7+zy^9$ &
$t + z^3 +zx+zy \subset \PP(3,1,2,2)$  & 28\\
& & $\PP^2$,\
$\Diff$=\scriptsize{$(0,0,\frac67,\frac89)$}+
$\frac12\GG+\frac34\GGG$ & \\
\hline
13& $zx^7+zxy^7$ &
$t + z^3 +zx^7+zxy \subset \PP(21,7,2,12)$  & 8\\
& & $\PP(7,1,6)$,\
$\Diff$=\scriptsize{$(0,0,0,\frac67)$}+
$\frac12\GG+\frac34\GGG$ & \\
\hline
14& $zx^7+zx^3y^5$ &
$t + z^3 +zx^7+zx^3y \subset \PP(21,7,2,8)$  & 16\\
& & $\PP(7,1,4)$,\
$\Diff$=\scriptsize{$(0,0,0,\frac45)$}+
$\frac12\GG+\frac34\GGG$ & \\
\hline
15& $zx^6y+azx^3y^5+bzy^9$ &
$t + z^3 +zx^2y+azxy^5+bzy^9 \subset \PP(27,9,8,2)$  & 12\\
& & $\PP(9,4,1)$,\
$\Diff$=\scriptsize{$(0,0,\frac23,0)$}+
$\frac12\GG+\frac34\GGG$ &\\
\hline
16& $zx^6y+zxy^n$ &
$t + z^3 +g(z,x,y) \subset \PP(18n-3,6n-1,2n-2,10)$  & 8,20\\
&$n=7;8$ & $\PP(6n-1,n-1,5)$,\
$\Diff$=
$\frac12\GG+\frac34\GGG$ & \\
\hline
17& $zx^5y^2+zy^n$ &
$t + z^3 +zxy^2+zy^n \subset \PP(3n,n,2n-4,2)$  & 16,20\\
&$n=9;11$ & $\PP(n,n-2,1)$,\
$\Diff$=\scriptsize{$(0,0,\frac45,0)$}+
$\frac12\GG+\frac34\GGG$ & \\
\hline
18& $zx^5y^2+zxy^n$ &
$t + z^3 +g(z,x,y) \subset \PP(15n-6,5n-2,2n-4,8)$  & 12,16\\
&$n=7;9$ & $\PP(5n-2,n-2,4)$,\
$\Diff$=
$\frac12\GG+\frac34\GGG$ & \\
\hline
19& $zx^5y^2+azx^3y^5+bzxy^8$ &
$t + z^3 +g(z,x,y) \subset \PP(57,19,6,4)$  & 8\\
& & $\PP(19,3,2)$,\
$\Diff$=
$\frac12\GG+\frac34\GGG$ &\\
\hline
20& $zx^5y^2+zx^2y^7$ &
$t + z^3 +g(z,x,y) \subset \PP(93,31,10,6)$  & 12\\
& & $\PP(31,5,3)$,\
$\Diff$=
$\frac12\GG+\frac34\GGG$ &\\
\hline
\end{supertabular}
}
\end{center}


\vspace{0.3cm}
\hoffset=-1cm

\begin{center}
\small{2. Singularity -- $t^2+z^4+g(z,x,y)$}
\end{center}

\vspace{0.5cm}

\begin{center}

\tabletail{\hline}
\tablehead{\hline}
\footnotesize{
\begin{supertabular}{|c|c|c|c|}
\hline
1& $zx^if_{4-i}(x,y^2)$ &
$t^2 + z^4 +zx^if_{4-i}(x,y) \subset \PP(8,4,3,3)$  & 3\\
& &
$\Diff$=\scriptsize{$(0,0,0,\frac12)$}+
$\frac23\GGG$ & \\
\hline
2& $zx^4+zy^n$ &
$t^2 + z^4 +zx+zy \subset \PP(2,1,3,3)$  & 5,8,12\\
&$n=5;7;11$ &
$\Diff$=\scriptsize{$(0,0,\frac34,\frac{n-1}{n})$}+
$\frac23\GGG$ &\\
\hline
3& $zx^4+azx^2y^5+bzy^{10}$ &
$t^2 + z^4 +zx^2+azxy+bzy^2 \subset \PP(4,2,3,3)$  & 6\\
& &
$\Diff$=\scriptsize{$(0,0,\frac12,\frac45)$}+
$\frac23\GGG$ &\\
\hline
4& $zx^4+zxy^n$ &
$t^2 + z^4 +zx^4+zxy \subset \PP(8,4,3,9)$  & 5,9\\
&$n=4,5;7,8$ &
$\Diff$=\scriptsize{$(0,0,0,\frac{n-1}{n})$}+
$\frac23\GGG$ &\\
\hline
5& $zx^if_{5-i}(x,y)$ &
$t^2 + z^4 +zx^if_{5-i}(x,y) \subset \PP(10,5,3,3)$  & 3\\
&$i=1$ &
$\Diff$=
$\frac23\GGG$ & \\
\hline
6& $zx^5+zy^7$ &
$t^2 + z^4 +zx+zy \subset \PP(2,1,3,3)$  & 15\\
& &
$\Diff$=\scriptsize{$(0,0,\frac45,\frac67)$}+
$\frac23\GGG$ &\\
\hline
7& $zx^5+zxy^5$ &
$t^2 + z^4 +zx^5+zxy \subset \PP(10,5,3,12)$  & 6\\
& &
$\Diff$=\scriptsize{$(0,0,0,\frac45)$}+
$\frac23\GGG$ &\\
\hline
8& $zx^5+zx^2y^4$ &
$t^2 + z^4 +zx^5+zx^2y \subset \PP(10,5,3,9)$  & 9\\
& &
$\Diff$=\scriptsize{$(0,0,0,\frac34)$}+
$\frac23\GGG$ &\\
\hline
9& $zx^4y+azx^2y^4+bzy^7$ &
$t^2 + z^4 +zx^2y+azxy^4+bzy^7 \subset \PP(14,7,9,3)$  & 6\\
& &
$\Diff$=\scriptsize{$(0,0,\frac12,0)$}+
$\frac23\GGG$ &\\
\hline
10& $zx^4y+zy^8$ &
$t^2 + z^4 +zxy+zy^8 \subset \PP(16,8,21,3)$  & 12\\
& &
$\Diff$=\scriptsize{$(0,0,\frac34,0)$}+
$\frac23\GGG$ &\\
\hline
11& $zx^4y+zxy^n$ &
$t^2 + z^4 +g(z,x,y) \subset \PP(8n-2,4n-1,3n-3,9)$  & 6,9\\
&$n=5;6$ &
$\Diff$=
$\frac23\GGG$ & \\
\hline

\end{supertabular}
}
\end{center}

\vspace{0.3cm}

\begin{center}
\small{3. Singularity -- $t^2+z^3x+g(z,x,y)$}
\end{center}

\begin{center}

\tabletail{\hline}
\tablehead{\hline}
\footnotesize{
\begin{supertabular}{|c|c|c|c|}
\hline
1& $x^5+xy^n$ &
$t + zx +x^5+xy \subset \PP(5,4,1,4)$  & 8,12,16\\
&$n=7;11;13;17;19,23$ & $\PP^2$,\
$\Diff$=\scriptsize{$(0,\frac23,0,\frac{n-1}{n})$}+
$\frac12\GG+\frac78\GGGG$ & 18,24\\
\hline
2& $x^5+azx^3y^3+bzxy^{12}+$ &
$t + z^3x +g(z,x,y^{1/3}) \subset \PP(15,4,3,2)$  & 4\\
&$cx^3y^9+dxy^{18}, |b|+|d|\ne 0$ & $\PP(2,3,1)$,\
$\Diff$=\scriptsize{$(0,0,0,\frac23)$}+
$\frac12\GG+\frac34\GGGG$ &\\
\hline
3& $x^5+azxy^{2n}+bxy^{3n}$ &
$t + z^3x +g(z,x,y^{1/n}) \subset \PP(15,4,3,4)$  & 8\\
&$n=3,5,7$ & $\PP(1,3,1)$,\
$\Diff$=\scriptsize{$(0,0,0,\frac{n-1}{n})$}+
$\frac12\GG+\frac78\GGGG$ &\\
\hline
4& $x^5+ax^3y^n+bxy^{2n}$ &
$t + zx +x^5+ax^3y+bxy^2 \subset \PP(5,4,1,2)$  & 6,12\\
&$n=5\ b\ne 0;7,11$ & $\PP(2,1,1)$,\
$\Diff$=\scriptsize{$(0,\frac23,0,\frac{n-1}{n})$}+
$\frac12\GG+\frac34\GGGG$ & \\
\end{supertabular}
}
\end{center}

\vspace{0.3cm}
\begin{center}
\small{4. Singularity -- $t^2+z^3y+g(z,x,y)$}
\end{center}

\begin{center}

\tabletail{\hline}
\tablehead{\hline}
\footnotesize{
\begin{supertabular}{|c|c|c|c|}
\hline
1& $x^7y+y^n$ &
$t + zy +xy+y^n \subset \PP(n,n-1,n-1,1)$  & 16,21,36\\
&$n=9;11;17;21$ & $\PP^2$,\
$\Diff$=\scriptsize{$(0,\frac23,\frac67,0)$}+
$\frac12\GG+\frac{2n-3}{2n-2}\GGGGG$ & 42\\
\hline
2& $x^7y+azy^{4n+3}+by^{6n+4}$ &
$t^2 + z^3y +g(z,x^{1/7},y) \subset \PP(3n+2,2n+1,6n+3,1)$  & 7\\
&$n=1,2$ &
$\Diff$=\scriptsize{$(0,0,\frac67,0)$}+
$\frac{2n}{2n+1}\GGGGG$ &\\
\hline
3& $x^7y+y^{2n}$ &
$t^2 + zy + xy+y^{2n}\subset \PP(n,2n-1,2n-1,1)$  & 15,21\\
&$n=6,7;9,10$ &
$\Diff$=\scriptsize{$(0,\frac23,\frac67,0)$}+
$\frac{2n-2}{2n-1}\GGGGG$ &\\
\hline
4& $x^7y+azy^{2n+1}+by^{3n+1}$ &
$t + z^3y +g(z,x^{1/7},y) \subset \PP(3n+1,n,3n,1)$  & 8,14\\
&$n=4;6$ & $\PP(1,3,1)$,\
$\Diff$=\scriptsize{$(0,0,\frac67,0)$}+
$\frac12\GG+\frac{2n-1}{2n}\GGGGG$ & \\
\end{supertabular}
}
\end{center}

\vspace{0.3cm}
\begin{center}
\small{5. Singularity -- $t^2+g_5+g(z,x,y)$.}
\end{center}

\begin{center}

\tabletail{\hline}
\tablehead{\hline}
\footnotesize{
\begin{supertabular}{|c|c|c|c|}
\hline
1& $z^5+zx^if_{5-i}(x,y)$ &
$t +z^5+ zx^if_{5-i}(x,y) \subset \PP(25,5,4,4)$  & 8\\
& & $\PP(5,1,1)$,\
$\Diff$=$\frac12\GG+\frac78\GGG$ &\\
\hline
2& $z^5+zx^3y^2 +az^3y^3+$ &
$t  +z^5+zxy^2 +az^3y^3+bzy^6 \subset \PP(15,3,8,2)$  & 12\\
&$bzy^6$ & $\PP(3,4,1)$,\
$\Diff$=\scriptsize{$(0,0,\frac23,0)$}+
$\frac12\GG+\frac34\GGG$ & \\
\hline
3& $z^5+zx^3y^2 +zy^7$ &
$t  +z^5+zxy^2 +zy^7 \subset \PP(35,7,20,4)$  & 24\\
& & $\PP(7,5,1)$,\
$\Diff$=\scriptsize{$(0,0,\frac23,0)$}+
$\frac12\GG+\frac78\GGG$ & \\
\hline
4& $z^5+zx^3y^2 +zxy^5$ &
$t  +z^5+zx^3y^2 +zxy^5 \subset \PP(65,13,12,8)$  & 16\\
& & $\PP(13,3,2)$,\
$\Diff$=$\frac12\GG+\frac78\GGG$ & \\
\end{supertabular}
}
\end{center}

\vspace{0.3cm}
\begin{center}
\small{6. Singularity -- $t^3+g(z,x,y)$.}
\end{center}

\begin{center}

\tabletail{\hline}
\tablehead{\hline}
\footnotesize{
\begin{supertabular}{|c|c|c|c|}
\hline
1& $z^2xy+x^5+xy^5$ &
$t + z^2xy+x^5+xy^5\subset \PP(25,8,5,4)$  & 12\\
& & $\PP(2,5,1)$,\
$\Diff$=$\frac23\GG+\frac{11}{12}\GGGG$ &\\
\hline
2& $z^2xy+x^4y +aty^4+$ &
$t^3+zxy+x^4y +aty^4+by^6 \subset \PP(8,15,5,4)$  & 8\\
&$+by^6$ &
$\Diff$=\scriptsize{$(0,\frac12,0,0)$}+
$\frac45\GGGGG$ & \\
\hline
3& $z^2xy+x^4y +ax^2y^4+$ &
$t  +zxy+x^4y +ax^2y^4+by^7 \subset \PP(14,9,3,2)$  & 12\\
&$+by^7$ &$\PP(3,1,2)$,\
$\Diff$=\scriptsize{$(0,\frac12,0,0)$}+
$\frac23\GG+\frac89\GGGGG$ & \\
\hline
4& $z^2xy+x^4y +y^8$ &
$t  +zxy+x^4y +y^8 \subset \PP(32,21,7,4)$  & 24\\
& &$\PP(3,1,4)$,\
$\Diff$=\scriptsize{$(0,\frac12,0,0)$}+
$\frac23\GG+\frac{20}{21}\GGGGG$ & \\
\end{supertabular}
}
\end{center}

\vspace{0.3cm}
\begin{center}
\small{7. Singularity -- $t^2z+g(z,x,y)$.}
\end{center}

\begin{center}

\tabletail{\hline}
\tablehead{\hline}
\footnotesize{
\begin{supertabular}{|c|c|c|c|}
\hline
1& $z^ix^jf_{4-i-j}(z,x)+$ &
$t^2z +g(z,x,y^{1/2}) \subset \PP(3,2,2,3)$  & 3\\
&$+axy^4+bzy^4+ctxy^2$ &
$\Diff$=\scriptsize{$(0,0,0,\frac12)$}+
$\frac23\ZZZZ$ &\\
&$(i\le 1 \& c\ne 0)\&(j\le 1 || b\ne 0)$&&\\
\hline
2& $z^ix^jf_{4-i-j}(z,x)+$ &
$tz +g(z,x,y^{1/5}) \subset \PP(3,1,1,3)$  & 6\\
&$+(zy^5|||xy^5)$ &
$\Diff$=\scriptsize{$(\frac12,0,0,\frac45)$}+
$\frac23\ZZZZ$ &\\
&$(j\le 2)\&(i=0 ||| i\le 1)$&&\\
\hline
3& $z^4+yx^if_{3-i}(x,y^3)+$ &
$tz +g(z,x,y) \subset \PP(15,5,6,2)$  & 4\\
&$+(az^2xy^2|||az^2y^5)$ &
$\Diff$=\scriptsize{$(\frac12,0,0,0)$}+
$\frac23\ZZZZZ$ &\\
\end{supertabular}
}
\end{center}

\vspace{0.3cm}
\begin{center}
\small{8. Singularity -- $t^2x+g(z,x,y)$.}
\end{center}

\begin{center}

\tabletail{\hline}
\tablehead{\hline}
\footnotesize{
\begin{supertabular}{|c|c|c|c|}
\hline
1& $z^4+x^3y+xy^6$ &
$t^2x +z+x^3y+xy^6 \subset \PP(6,17,5,2)$  & 8\\
& &$\PP(3,5,1)$,\
$\Diff$=$\frac34\ZZZ+\frac78\ZZZZ$ &\\
\hline
2& $z^4+x^5+xy^5$ &
$t^2x +z+x^5+xy \subset \PP(2,5,1,4)$  & 16\\
& &$\PP(1,1,2)$,\
$\Diff$=\scriptsize{$(0,0,0,\frac45)$}+
$\frac34\ZZZ+\frac78\ZZZZ$ &\\
\hline
3& $z^4+x^4y+xy^5$ &
$t^2x +z^2+x^4y+xy^5 \subset \PP(15,19,8,6)$  & 12\\
& &
$\Diff$=$\frac12\ZZZ+\frac56\ZZZZ$ &\\
\hline
4& $z^4+x^4y+zy^n$ &
$tx +g(z,x,y) \subset \PP(12n+3,4n,4n-3,12)$  & 10,24\\
&$n=4;5$ &
$\Diff$=\scriptsize{$(\frac12,0,0,0)$}+
$\frac23\ZZZZ$ &\\
\end{supertabular}
}
\end{center}

\vspace{0.3cm}
\begin{center}
\small{9. Singularity -- $t^2y+g(z,x,y)$.}
\end{center}

\begin{center}

\tabletail{\hline}
\tablehead{\hline}
\footnotesize{
\begin{supertabular}{|c|c|c|c|}
\hline
1& $z^4+zx^4+az^2xy^2+$ &
$ty +g(z,x,y) \subset \PP(27,8,6,5)$  & 10\\
&$+bx^2y^4$ &
$\Diff$=\scriptsize{$(\frac12,0,0,0)$}+
$\frac23\ZZZZZ$ &\\
\hline
2& $z^4+zx^4+zy^5$ &
$t^2y +z^4+zx^2+zy^5 \subset \PP(17,10,15,6)$  & 12\\
& &
$\Diff$=\scriptsize{$(0,0,\frac12,0)$}+
$\frac23\GGG$ &\\
\hline
3& $z^4+zx^4+zxy^4$ &
$ty +z^4+zx^4+zxy^4 \subset \PP(55,16,12,9)$  & 18\\
& &
$\Diff$=$\frac12\GG+\frac56\GGG$ &\\
\hline
4& $z^3x+x^5+azx^2y^3+$ &
$t^2y +g(z,x,y) \subset \PP(20,12,9,5)$  & 5\\
&$+by^9$ &
$\Diff$=$\frac34\XXX$ &\\
\hline
5& $z^3x+x^5+y^{13}$ &
$t^2y +zx+x^5+y^{13} \subset \PP(30,52,13,5)$  & 15\\
& &
$\Diff$=\scriptsize{$(0,\frac23,0,0)$}+
$\frac12\XXX$ &\\
\hline
6& $z^3x+x^5+xy^n$ &
$ty +zx+x^5+xy^n \subset \PP(5n-4,4n,n,4)$  & 9,16,24\\
&$n=5;7;11$ &
$\Diff$=\scriptsize{$(0,\frac23,0,0)$}+
$\frac12\GG+\frac78\GGGG$ &\\
\hline
7& $z^3x+x^5+azx^3y+$ &
$ty +g(z,x,y) \subset \PP(13,4,3,2)$  & 4\\
&$bzxy^4+cx^3y^3+dxy^6$ &
$\Diff$=$\frac12\GG+\frac34\GGGG$ &\\
&$|b|+|d|\ne 0$&&\\
\hline
8& $z^3x+x^5+azxy^6+$ &
$ty +g(z,x,y) \subset \PP(41,12,9,4)$  & 8\\
&$+bxy^9$ &
$\Diff$=$\frac12\GG+\frac78\GGGG$ &\\
\hline
9& $z^3x+x^5+ax^3y^5+$ &
$ty +g(z^{1/3},x,y) \subset \PP(23,20,5,2)$  & 12\\
&$+bxy^{10}$ &
$\Diff$=\scriptsize{$(0,\frac23,0,0)$}+
$\frac12\GG+\frac34\GGGG$ &\\
\hline
10& $z^3x+x^5+x^2y^7$ &
$t^2y + zx+x^5+x^2y^7\subset \PP(16,28,7,3)$  & 9\\
& &
$\Diff$=\scriptsize{$(0,\frac23,0,0)$}+
$\frac34\XXX$ &\\
\hline
11& $z^3x+x^5+zy^{2n+1}$ &
$t^2y + g(z,x,y)\subset \PP(15n+2,8n+4,6n+3,11)$  & 5,11\\
&$n=2;4^*$ &
$\Diff$=$\frac34\XXX$; $(\frac12\XXX)^*$ &\\
\hline
12& $z^3x+x^5+zxy^n$ &
$ty +g(z,x,y) \subset \PP(15n-8,4n,3n,8)$  & 8,10,16\\
&$n=3;5;7$ &
$\Diff$=$\frac12\GG+\frac78\GGGG$ &\\
\hline
13& $z^3x+x^5+z^2y^5$ &
$t^2y + z^3x+x^5+z^2y^5\subset \PP(34,20,15,7)$  & 7\\
& &
$\Diff$=$\frac12\XXX$ &\\
\hline
14& $z^3x+x^6+azxy^4+$ &
$ty +g(z,x,y) \subset \PP(31,10,6,5)$  & 10\\
&$+bxy^6$ &
$\Diff$=$\frac12\GG+\frac9{10}\GGGG$ &\\
\hline
15& $z^3x+x^6+xy^7$ &
$ty + zx+x^6+xy^7\subset \PP(37,35,7,5)$  & 30\\
& &
$\Diff$=\scriptsize{$(0,\frac23,0,0)$}+
$\frac12\GG+\frac9{10}\GGGG$ &\\
\hline
16& $z^3x+x^5y+y^{2n+1}$ &
$t^2y +g(z^{1/3},x,y)\subset \PP(5n,8n+5,2n,5)$  & 9,15\\
&$n=3;4$ &
$\Diff$=$\frac23\ZZZ+\frac{3n-1}{3n}\ZZZZZ$ &\\
\hline
17& $z^3x+x^5y+azx^2y^3+$ &
$ty +g(z,x,y)\subset \PP(35,11,7,5)$  & 10\\
&$+by^8$ &
$\Diff$=\scriptsize{$(\frac12,0,0,0)$}+
$\frac67\ZZZZZ$ &\\
\hline
18& $z^3x+x^5y+y^{10}$ &
$ty +zx+x^5y+y^{10}\subset \PP(45,41,9,5)$  & 30\\
& &
$\Diff$=\scriptsize{$(\frac12,0,0,0)$}+
$\frac23\ZZZ+\frac{26}{27}\ZZZZZ$ &\\
\hline
19& $z^3x+x^5y+azxy^4+$ &
$ty +g(z,x,y)\subset \PP(25,8,5,4)$  & 8\\
&$+bxy^6$ &
$\Diff$=$\frac12\GG+\frac78\GGGG+\frac45\ZZZZZ$ &\\
\hline
20& $z^3x+x^5y+ax^3y^4+$ &
$ty +g(z^{1/3},x,y)\subset \PP(15,14,3,2)$  & 12\\
&$+bxy^7$ &
$\Diff$=$\frac12\GG+\frac23\ZZZ+\frac34\GGGG+\frac89\ZZZZZ$ &\\
\hline
21& $z^3x+x^5y+xy^8$ &
$ty +zx+x^5y+xy^8\subset \PP(35,32,7,4)$  & 24\\
& &
$\Diff$=$\frac12\GG+\frac23\ZZZ+\frac78\GGGG+\frac{20}{21}\ZZZZZ$ &\\
\hline
22& $z^3x+x^5y+x^2y^5$ &
$t^2y +zx+x^5y+x^2y^5\subset \PP(10,19,4,3)$  & 9\\
& &
$\Diff$=$\frac23\ZZZ+\frac56\ZZZZZ$ &\\
\hline
23& $z^3x+x^5y+x^2y^6$ &
$ty +zx+x^5y+x^2y^6\subset \PP(25,23,5,3)$  & 18\\
& &
$\Diff$=\scriptsize{$(\frac12,0,0,0)$}+
$\frac23\ZZZ+\frac{14}{15}\ZZZZZ$ &\\
\hline
24& $z^3x+x^5y+zy^n$ &
$ty +g(z,x,y)\subset \PP(15n-10,4n+1,3n-2,11)$  & 14,22\\
&$n=5;7$ &
$\Diff$=\scriptsize{$(\frac12,0,0,0)$}+
$\frac{3n-3}{3n-2}\ZZZZZ$ &\\
\hline
25& $z^3x+x^5y+zy^6$ &
$t^2y +z^3x+x^5y+zy^6\subset \PP(40,25,16,11)$  & 11\\
& &
$\Diff$=$\frac45\XXX+\frac78\ZZZZZ$ &\\
\hline
26& $z^3x+x^5y+zxy^5$ &
$ty +z^3x+x^5y+zxy^5\subset \PP(65,20,13,8)$  & 16\\
& &
$\Diff$=$\frac12\GG+\frac78\ZZZZ+\frac45\XXX$ &\\
\hline
27& $z^3x+x^5y+z^2y^3$ &
$t^2y +z^3x+x^5y+z^2y^3\subset \PP(20,13,8,7)$  & 7\\
& &
$\Diff$=$\frac34\ZZZZZ$ &\\
\hline
28& $z^3x+x^5y+z^2y^4$ &
$ty +z^3x+x^5y+z^2y^4\subset \PP(55,17,11,7)$  & 14\\
& &
$\Diff$=\scriptsize{$(\frac12,0,0,0)$}+
$\frac{10}{11}\ZZZZZ$ &\\
\hline
27& $z^3x+zx^4+y^7$ &
$t^2y +z^3x+zx^4+y^7\subset \PP(33,21,14,11)$  & 7\\
& &
$\Diff$=$\frac23\XXX$ &\\
\hline
28& $z^3x+zx^4+y^9$ &
$t^2y +z^3x+zx^4+y^9\subset \PP(44,27,18,11)$  & 11\\
& &
$\Diff$=$\frac12\ZZZZZ$ &\\
\hline
29& $z^3x+zx^4+y^{10}$ &
$ty +z^3x+zx^4+y^{10}\subset \PP(99,30,20,11)$  & 22\\
& &
$\Diff$=\scriptsize{$(\frac12,0,0,0)$}+
$\frac23\XXX$ &\\
\hline
30& $z^3x+zx^4+xy^{2n+1}$ &
$t^2y +g(z,x,y)\subset \PP(11n+1,6n+3,4n+2,9)$  & 5,9\\
&$n=2;3$ &
$\Diff$=$\frac23\GGGG+\frac{n-2}{n-1}\ZZZZZ$ &\\
\hline
31& $z^3x+zx^4+azxy^4+$ &
$ty +g(z,x,y)\subset \PP(19,6,4,3)$  & 6\\
&$+bx^4y^2+cxy^6$ &
$\Diff$=$\frac12\GG+\frac56\GGGG$ &\\
\hline
32& $z^3x+zx^4+xy^8$ &
$ty +z^3x+zx^4+xy^8\subset \PP(79,24,16,9)$  & 18\\
& &
$\Diff$=$\frac12\GG+\frac56\GGGG$ &\\
\hline
33& $z^3x+zx^4+x^2y^5$ &
$t^2y +z^3x+zx^4+x^2y^5\subset \PP(24,15,10,7)$  & 7\\
& &
$\Diff$=$\frac12\ZZZZZ$ &\\
\hline
34& $z^3x+zx^4+az^2y^3+$ &
$t^2y +g(z,x,y)\subset \PP(14,9,6,5)$  & 5\\
&$+bx^3y^3, a\ne 0$ &
$\Diff$=$\frac12\ZZZZZ$ &\\
\hline
35& $z^3x+zx^4+az^2y^4+$ &
$ty +g(z,x,y)\subset \PP(39,12,8,5)$  & 10\\
&$+bx^3y^4$ &
$\Diff$=\scriptsize{$(\frac12,0,0,0)$}+
$\frac23\XXX$ &\\
\hline
36& $z^3x+zx^4+zy^n$ &
$ty +g(z,x,y)\subset \PP(11n-8,3n,2n,8)$  & 10,16\\
&$n=5;7$ &
$\Diff$=$\frac12\GG+\frac34\GGG+\frac{n-5}{n-4}\XXX$ &\\
\hline
37& $z^3x+zx^4+azx^2y^3+$ &
$ty +g(z,x,y)\subset \PP(29,9,6,4)$  & 8\\
&$+bzy^6$ &
$\Diff$=$\frac12\GG+\frac34\GGG$ &\\
\hline
38& $z^3x+zx^4+zxy^5$ &
$ty +z^3x+zx^4+zxy^5\subset \PP(49,15,10,6)$  & 12\\
& &
$\Diff$=$\frac12\GG+\frac34\GGG+\frac56\GGGG$ &\\
\hline
39& $z^3x+tx^3+ax^3y^4+$ &
$t^2y +g(z^{1/3},x,y)\subset \PP(12,22,5,3)$  & 9\\
&$+by^9$ &
$\Diff$=\scriptsize{$(0,\frac23,0,0)$}+
$\frac12\XXX$ &\\
\hline
40& $z^3x+tx^3+xy^n$ &
$t^2y +g(z^{1/3},x,y)\subset \PP(3n-2,5n,n+1,5)$  & 9,12\\
&$n=5;7$ &
$\Diff$=\scriptsize{$(0,\frac23,0,0)$}+
$\frac45\GGGG$ &\\
\hline
41& $z^3x+tx^3+azxy^4+$ &
$t^2y +g(z,x,y)\subset \PP(16,10,7,5)$  & 5\\
&$+bxy^6$ &
$\Diff$=$\frac45\GGGG+\frac12\XXX$ &\\
\hline
42& $z^3x+tx^3+xy^8$ &
$t^2y +zx+tx^3+xy^8\subset \PP(22,40,9,5)$  & 15\\
& &
$\Diff$=\scriptsize{$(0,\frac23,0,0)$}+
$\frac45\GGGG+\frac12\XXX$ &\\
\hline
43& $z^3x+tx^3+zy^{2n+1}$ &
$t^2y +g(z,x,y)\subset \PP(18n+2,10n+4,6n+5,13)$  & 8,13\\
&$n=2;3$ &
$\Diff$=$\frac12\XXX$ &\\
\hline
44& $z^3x+tx^3+zxy^n$ &
$t^2y +g(z,x,y)\subset \PP(9n-4,5n,3n+2,10)$  & 5,10\\
&$n=3;5$ &
$\Diff$=$\frac45\GGGG$ &\\
\hline
45& $z^3x+tx^3+zx^2y^3$ &
$t^2y +z^3x+tx^3+zx^2y^3\subset \PP(26,16,11,7)$  & 7\\
& &
$\Diff$=$\frac12\XXX$ &\\
\hline

\end{supertabular}
}
\end{center}

\end{document}